\documentclass[12pt]{article}
\usepackage{subeqn}

\newcommand{\eqref}[1]{\textup{(\ref{#1})}}

\def\Boxx{
\vbox{
\halign to5.8pt %5pt
{\strut##&
\hfil ## \hfil \cr
&$\kern -0.5pt
\sqcap$ \cr
\noalign{\kern -5pt%-7.5pt
\hrule}
}}}

\def\bea{\begin{eqnarray}}
\def\eqnn{\bea\label}
\def\eea{\end{eqnarray}}
\def\nn{\nonumber}

\newcommand{\eqna}[1]{\begin{subequations} \label{#1}
\begin{eqnarray}}
\def\eena{\end{eqnarray}
\end{subequations}}

\def\hF{{\hat F}} \def\tF{{\tilde F}}
\def\hJ{{\hat J}} \def\tJ{{\tilde J}} \def\tI{{\tilde I}}
\def\hK{{\hat K}} \def\tK{{\tilde K}}
\def\hp{{\hat p}}\def\hr{{\hat r}}
\def\tp{{\tilde p}} \def\tr{{\tilde r}}
 \def\y{\eta}
 \def\s{\sigma}  \def\t{\tau}
\def\eps{\epsilon} \def\ve{\varepsilon}
\def\ra{\rightarrow}
\def\hvf{\hat{\varphi}}  
\def\tvf{\tilde{\varphi}}  
\def\bbz{Z\!\!\!Z}
\def\k{\kappa}
\def\Uf{U_q(sl(4))}
\def\hd{{\hat{\cal D}}}
\def\({\left(}
\def\){\right)}
\def\bbn{I\!\!N}
\def\bz{{\bar z}}
\def\G{\Gamma}
\def\b{\beta}
\def\cl{{\cal L}}
\def\be{\begin{equation}}
\def\eqn{\begin{equation}\label}
\def\ee{\end{equation}}

\def\bea{\begin{eqnarray}}
\def\eqnn{\bea\label}
\def\eea{\end{eqnarray}}
\def\nn{\nonumber}

\def\nt{\noindent}

\def\half{{\textstyle{1\over2}}}
\def\ha{{\textstyle{1\over2}}}
\def\sixth{{\textstyle{1\over6}}}
\def\bv{{\bar v}}
\def\hm{{\hat M}}
\def\tg{\tilde \gamma}
\def\hg{\hat \gamma}
\def\vf{\varphi} \def\td{{\tilde d}}
\def\pd{\partial} \def\tb{{\tilde\beta}}
\def\om{\omega} \def\th{{\tilde h}} \def\hh{{\hat h}}
\def\l{\lambda}

\begin{document}
\begin{center}
\vspace*{1.0cm}

{\LARGE{\bf $q$-Conformal Invariant Equations
and\\[2mm] $q$-Plane Wave Solutions}}

\vskip 1.5cm

{\large {\bf V.K. Dobrev}$^{1,2}$ ~and~ {\bf S.T. Petrov}$^{2,3}$}

\vskip 0.5 cm

$^1$ School of Informatics,\\
University of Northumbria,\\ Newcastle upon Tyne NE1 8ST, UK

$^2$ Institute of Nuclear Research and Nuclear Energy,\\
Bulgarian Academy of Sciences,\\ 72 Tsarigradsko Chaussee, 1784
Sofia, Bulgaria\\
(permanent address of both authors)

$^3$ Laboratoire de Physique des Mat\'eriaux,\\
Universit\'e  Henri Poincar\'e Nancy-I\\
B.P. 239, 54506 Vandoeuvre les Nancy,
Nancy, France

\end{center}

\vspace{1 cm}

\begin{abstract}
We give new solutions of the quantum conformal deformations of the
full Maxwell equations in terms of deformations of the plane wave.
We study the compatibility of these solutions with the
conservation of the current.
We also start the study of quantum linear conformal
(Weyl) gravity by writing the corresponding q-deformed equations.
\end{abstract}

\vspace{1 cm}

\section*{Introduction}

One of the purposes of quantum deformations
is to provide an alternative of the regularization
procedures of quantum field theory. Applied to
Minkowski space-time the quantum deformations approach
is also an alternative to Connes' noncommutative geometry \cite{Con}.
The first step in such an approach is to construct
a noncommutative quantum deformation of Minkowski space-time.
There are several possible such deformations,
cf. \cite{CSSW,SWZ,Mja,Mjb,VKD2}.
We shall follow the deformation of \cite{VKD2}
which is different from the others, the most
important aspect being that it is related
to a deformation of the conformal group.

The first problem to tackle in a noncommutative
deformed setting is to study the q-deformed analogues of
the conformally invariant equations. Here we continue the
study of hierarchies of deformed equations
derived in \cite{VKD2,VKD3,VKD4} with the use of
quantum conformal symmetry.
We give now a description of our setting starting from the
simplest example.

It is well known that the d'Alembert equation
\eqn{dal} \Boxx ~\vf(x) ~~=~~ 0 ~, \qquad \Boxx ~=~ \pd^\mu \pd_\mu
~=~ (\pd_0)^2 - (\vec\pd)^2 ~, \end{equation}
is conformally invariant,
cf., e.g., \cite{BR}. Here ~$\vf$~ is a scalar field of fixed conformal
weight, ~$x\ =\ (x_0,x_1,x_2,x_3)$~ denotes the Minkowski
space-time coordinates. Not known was the fact that (\ref{dal})\
may be interpreted as conditionally conformally invariant
equation and thus may be rederived from a subsingular vector of a
Verma module of the algebra\ $sl(4)$, the complexification of the
conformal algebra\ $su(2,2)$ \cite{VKD3}.

The same idea was used in \cite{VKD3} to derive a q-d'Alembert
equation, namely, as arising from a subsingular vector of a
Verma module of the quantum algebra\ $\Uf$. The resulting
equation is a q-difference equation and the solution spaces are
built on the noncommutative q-Minkowski space-time of \cite{VKD2}.

Besides the q-d'Alembert equation in \cite{VKD3} were derived a
whole hierarchy of equations corresponding to the massless
representations of the conformal group and parametrized by a
nonnegative integer \ $r$\ \cite{VKD3}. The case $r=0$ corresponds
to the q-d'Alembert equation, while for each \ $r>0$ \ there are
two couples of equations involving fields of conjugated Lorentz
representations of dimension $r+1$. For instance, the case ~$r=1$~
corresponds to the massless Dirac equation, one couple of
equations describing the neutrino, the other couple of equations
describing the antineutrino, while the case ~$r=2$~ corresponds to
the Maxwell equations.

The construction of solutions of the q-d'Alembert
hierarchy was started in
 \cite{DK} with the q-d'Alembert equation. One of the solutions given
was a deformation of the plane wave as a formal power series in
the noncommutative coordinates of q-Minkowski space-time and
four-momenta. This q-plane wave has some properties analogous to
the classical one but is not an exponent or q-exponent.
Thus, it differs conceptually from the classical
plane wave and may serve as a regularization of the latter.
In the same sense it differs from the q-plane wave in the
paper \cite{Mey}, which is not surprising, since there is used different
q-Minkowski space-time (from \cite{CSSW,SWZ,Mja} and different
q-d'Alembert equation both based only on a (different) q-Lorentz
algebra, and not on q-conformal (or $\Uf$) symmetry as in our
case. In fact, it is not clear whether the q-Lorentz algebra of
\cite{CSSW,SWZ,Mja} used in \cite{Mey} is extendable to a
q-conformal algebra.

For
the equations labelled by $r>0$ it turned out that one needs a
second q-deformation of the plane wave in a conjugated basis
 \cite{DGPZ}. The solutions of the hierarchy in terms of the two
q-plane waves were given in \cite{DGPZ} for $r=1$ and in
 \cite{DZ} for $r>1$. Later these two q-plane waves were generalized
and correspondingly more general solutions of the hierarchy were
given in \cite{DPZa}.

Another hierarchy derived in \cite{VKD2} is the Maxwell hierarchy.
The two hierarchies have only one common member
- the Maxwell equations - they are the lowest member of the
Maxwell hierarchy and the $r=2$ member of the massless hierarchy.
The compatibility of the solutions of the free q-Maxwell equations
with the q-potential equations was studied \cite{DPZb}.

In the present paper we study the full q-Maxwell equations and
the compatibility of their solutions with the conservation of the
current. We give new solutions of the full q-Maxwell equations in
two conjugated bases. The solutions of the homogeneous equations
are also new (generalizing previously given solutions).

Another family contained in \cite{VKD4}, but not explicated there,
is related with the linear conformal (Weyl) gravity which we
start to study in this paper. Namely, in the last Section we
write down the quantum conformal deformations of the linear conformal
(Weyl) gravity.

\section{Preliminaries}
\setcounter{equation}{0}

\nt
First we introduce new Minkowski variables:
\eqn{mink} x_\pm \ \equiv \ x_0 \pm x_3\ , \qquad v \ \equiv \ x_1 -i
x_2\ , \qquad \bv \ \equiv \ x_1 + i x_2\ , \ee
which, (unlike the $x_\mu$), have definite group--theoretical
interpretation as part of a six-dimensional coset of the conformal group\
$SU(2,2)$ (as explained in \cite{VKD2}). In terms of these
variables, e.g., the d'Alembert equation (\ref{dal}) is:
\eqn{dali} \Boxx\ \vf = (\pd_- \pd_+ \ -\
\pd_{v}\ \pd_\bv )\ \vf \ =\ 0 \ .\ee

In the q-deformed case we use
the noncommutative q-Minkowski space-time of
 \cite{VKD2} which is
given by the following commutation relations
(with $\l \equiv q-q^{-1}$):
\eqn{coop} x_\pm v = q^{\pm 1} v x_\pm \,, \ \
 x_\pm \bv = q^{\pm 1} \bv x_\pm \,,\ \
x_+ x_- - x_- x_+ = \l v\bv \,, \ \ \bv v = v\bv \,,
\ee
with the deformation parameter being a phase:\ $\vert q\vert =1$.
Relations (\ref{coop}) are preserved by the anti-linear anti-involution
 \ $\om$ \ :
\eqn{cnjm} \om (x_\pm) \ = \ x_\pm \ , \quad \om (v) \ = \ \bv \ ,
\quad \om (q) \ = \ \bar q \ = \ q^{-1} \ ,
\quad (\om (\l) \ = \ -\l ) \ . \ee

The solution spaces consist of formal power series in the
q-Minkowski coordinates (which we give in two conjugate bases):
\bea
&&\vf \ = \ \sum_{j,n,\ell,m\in\bbz_+}\ \mu_{j n\ell m}\
\vf_{j n\ell m} \ , \qquad \vf_{j n\ell m}\ =\
\hvf_{j n\ell m},\ \tvf_{j n\ell m} \ ,\label{spc}\\
&&\hvf_{j n\ell m} \ = \ v^j\ x^n_-\ x^\ell_+\ \bv^m \ ,\label{spca}\\ [2mm]
&&\tvf_{j n\ell m} \ = \ \bv^m\ x^\ell_+\ x^n_-\ v^j \ = \
\om (\hvf_{j n\ell m}) \ . \label{spcb}\eea
The solution spaces (\ref{spc}) are representation spaces of the
quantum algebra $U_q(sl(4))$. For the latter we use the rational
basis of Jimbo \cite{MJ}.
The action of\ $\Uf$\ on\ $\hvf_{j n\ell m}$\ was given in
 \cite{VKD1}, and on \ $\tvf_{j n\ell m}$ in \cite{DGPZ}.
Because of the conjugation $\om$ we are actually working with
the conformal quantum algebra which is a deformation of $U(su(2,2))$.

 Further we suppose that\ $q$\ is not a nontrivial root of unity.

In order to write our q-deformed equations in compact form
it is necessary to introduce some additional operators.
We first define the operators:
\eqnn{opes}
\hm^\pm_{\k}\ \vf \ =& \ \sum_{j,n,\ell,m\in\bbz_+}\ \mu_{j n\ell m}\
\hm^\pm_{\k}\ \vf_{j n\ell m} \ , \qquad \k = \pm, v,\bv \ ,
\\ [2mm]
T^\pm_{\k}\ \vf \ =& \ \sum_{j,n,\ell,m\in\bbz_+}\ \mu_{j n\ell m}\
T^\pm_{\k}\ \vf_{j n\ell m} \ , \qquad \k = \pm, v,\bv \ , \eea
and \ $\hm^\pm_{+}\,$, \ $\hm^\pm_{-}\,$, \ $\hm^\pm_{v}\,$,
 \ $\hm^\pm_{\bv}\,$, resp., acts on \ $\vf_{j n\ell m}$ \ by changing
by \ $\pm1$ \ the value of \ $j,n,\ell,m$, resp., while
 \ $T^\pm_{+}\,$, \ $T^\pm_{-}\,$, \ $T^\pm_{v}\,$,
 \ $T^\pm_{\bv}\,$, resp., acts on \ $\vf_{j n\ell m}$ \
by multiplication by \ $q^{\pm j},q^{\pm n},q^{\pm \ell},q^{\pm
m}$, resp. We shall use also the 'logs' $N_\k$ such that ~$T_\k
~=~ q^{N_\k}$. Now we can define the q-difference operators:
\eqn{qdif} \hd_{\k} \ \vf \ = \ {1\over \l}\ \hm^{-1}_{\k}\
\left( T_{\k} - T^{-1}_{\k} \right) \ \vf \  = \ {1\over \l}\ \hm^{-1}_{\k}\
\left( q^{N_{\k}} - q^{-N_{\k}} \right) \ \vf \ . \ee
Note that when $q\ra 1$ then $\hd_{\k}\ra \pd_k\,$.
Using (\ref{opes}) and (\ref{qdif}) the q-d'Alembert equation
may be written as \cite{VKD3}, \cite{DGPZ}, respectively,
\eqn{qdal}
\left( q\ \hd_-\ \hd_+\ T_v\ T_\bv \ -\ \hd_v\ \hd_\bv \right)\
T_v\ T_-\ T_+\ T_\bv \ \hvf \ \ = \ \ 0 \,, \ee
\eqn{qqdal}
\left( \hd_-\ \hd_+\ -\ q\ \hd_v\ \hd_\bv \ T_v\ T_\bv \right)\
T_-\ T_+\ \ \tvf \ \ = \ \ 0 \,. \ee
Note that when $q\ra 1$ both equations (\ref{qdal}), (\ref{qqdal})
go to (\ref{dali}). Note that the operators in (\ref{opes}), (\ref{qdif}),
(\ref{qdal}), (\ref{qqdal}) for different variables commute,
i.e., we have passed to commuting variables. However, keeping
the normal ordering it is straightforward to pass back to
noncommuting variables.

Next we recall that the Maxwell's equations are part also of
Maxwell's hierarchy of equations. The
quantum conformal deformation of the equations of
the hierarchy are \cite{VKD2}:
\be\label{mxg}
_qI^+_n ~ _qF^+_n ~=~ _qJ^n ~, \qquad
_qI^-_n ~ _qF^-_n ~=~ _qJ^n \ee
where in the basis (\ref{spca}) the operators are:
\bea
_qI^+_n ~&=&~ \half\ \Biggl(
\Bigl( q\ \hd_v\ ~+~
\hm_\bz\ \hd_{+}\ (T_{-}\ T_v)^{-1} ~T_{\bv}
\Bigr)\ T_-\, [n+2-N_z]_q ~- \nn \\ &&~-
~ q^{-n-2}\ \Bigl(\hd_{-}\ T_- ~+~
q^{-1}~ \hm_\bz ~\hd_{\bv} ~-~ \nn \\ &&~-
\l\ \hm_v\ \hm_\bz \ \hd_{-}\ \hd_{+}\ T_{\bv} \Bigr)
\ T_-^{-1}\, \hd_z \Biggr)\ T_+\ T_v\ T_z\ T_\bz^{-1}
\ ,\label{intp}\\
_qI^-_n ~&=&~ \half\ \Bigl( \hd_{\bv}
~+~ q\ \hm_z \ \hd_{+}\ T_{\bv}\ T_{-}\ T_v^{-1} ~-\nn\\ &&-~
q \l\ \hm_v\ \hd_{-}\ \hd_{+}\ T_{\bv}
\Bigr)\ T_{\bv}\ [n+2-N_\bz]_q
 ~-\nn\\ &&-~ \half\ q^{n+3} ~
\Bigl(\hd_{-} ~+~ q\ \hm_z\ \hd_v ~T_{-}
\Bigr)\ \hd_\bz \ T_{-}\ T_{\bv} \ , \label{intm} \eea
while where in the basis (\ref{spcb}) the operators are:
\bea
_qI^+_n ~&=&~ \half\
q\, \Bigl( \hd_v\ ~+~
\hm_\bz\ \hd_{+}\ T_{-}\ T_\bv^{-1} ~T_{v}
\Bigr)\ T_v \ [n+2-N_z]_q ~-\nn\\ &&-
~\half\ q^{n+3}\ \Bigl(\hd_{-} ~+~
 \hm_\bz ~\hd_{\bv}~T_- ~+\nn \\ &&~+~
\l\ q^{-1}\hm_v\ \hm_\bz \ \hd_{-}\ \hd_{+}\ T_{\bv}^{-1}\ T_- \Bigr)
\ \hd_z\ T_- \ T_v
\ ,\label{intps}\\
_qI^-_n ~&=&~ \half\ \Biggl(
\Bigl( \hd_{\bv}\ T_{\bv}\ T_{-}
~+~ \hm_z \ \hd_{+}\ T_{v} ~+\nn\\ &&+~
q^{-1} \l\ \hm_v\ \hd_{-}\ \hd_{+}\ T_{-}
\Bigr)\ [n+2-N_\bz]_q
 ~-\nn\\ &&-~ q^{-n-2} ~
\Bigl(\hd_{-} ~+~ \hm_z\ \hd_v ~T_{-}^{-1}
\Bigr)\ \hd_\bz\ T_\bv \Biggr)\ T_{+}\ T_{\bz}\ T_z^{-1} \ .
\label{intms} \eea
Note that for $q=1$ (\ref{intp}),(\ref{intm}) coincide
with (\ref{intps}),(\ref{intms}), respectively.
Maxwell's equations $\ \pd^\mu F_{\mu\nu} ~=~ J_\nu$,
$\ \eps_{\mu\nu\rho\s} \pd^\mu F^{\rho\s} \ =\ 0\ $
are obtained from (\ref{mxg}) for $n=0$,
$q=1$, substituting the fixed helicity constituents $\ F^\pm\
$ by: $\ F^+ ~=~ z^2 (F_1^+ +iF_2^+) -2z F_3^+ -
(F_1^+ - iF^+_2)$, $\ F^- ~=~ \bz^2 (F_1^- -iF_2^-) - 2\bz F_3^- -
(F_1^- +iF^-_2)$, $\ F^\pm_k = F_{k0} \pm {i\over 2} \ve_{k\ell m}
F_{\ell m} = E_k \pm i H_k\,$, $\ J^0 ~=~ \bz z (J_0 + J_3) + z (J_1 +i J_2)
+ \bz (J_1 -i J_2) + (J_0 - J_3)$, and then comparing the
coefficients of the resulting first order polynomials in $z$ and
$\bz$.

We shall look for solutions of the full q-Maxwell's equations in
terms of deformations of the plane wave.
Let us first recall these deformations from
 \cite{DPZa}. The first deformation is given in the basis (\ref{spca}):
\bea\label{dploska} \widehat{\exp}_q (k, x) \ &=&\
\sum_{s=0}^\infty \, {1 \over {[s]_q!}}\ \hh_s \ ,\\
&& [s]_q! \equiv [s]_q [s-1]_q \cdots [1]_q \ , \quad [0]_q! \equiv 1 \ ,
\quad [n]_q \equiv { q^n-q^{-n}\over q-q^{-1} } \ ,
\nn\eea
\bea\label{hfs}
\hh_s\ &=&\ {\footnotesize
 \beta^s \sum_{a,b,n \in \bbz_+}
{ (-1)^{s-a-b}\ q^{n(s-2a-2b+2n)\ +\ a(s-a-1)\ +\
b(-s+a+b+1)}\ q^{P_s(a,b)}
\over {\G_q(a-n+1)\G_q(b-n+1)\G_q(s-a-b+n+1)[n]_q! }} \ \times }\nn\\
&&\times \ k_v^{s-a-b+n} k_-^{b-n} k_+^{a-n} k_\bv^n v^n x_-
^{a-n} x_+^{b-n} \bv^{s-a-b+n} \ , \\[2mm]
&&\left(\beta^s \right)^{-1} \ = \ \sum_{p=0}^s \
{q^{(s-p)(p-1)+p} \over [p]_q!\ [s-p]_q!} \ ,\nn \eea
where the momentum components \ $(k_v,k_-,k_+,k_\bv)$\
are supposed to be non-commu\-ta\-tive between themselves (obeying the same
rules (\ref{coop}) as the q-Min\-kow\-ski coordinates), and commutative
with the coordinates. Further, $\G_q$ is a $q$-deformation of the
$\G$-function, of which here we use only the properties:
$\G_q(p) = [p-1]_q!$ for $p\in\bbn$,
$1/\G_q(p) = 0$ for $p\in\bbz_-\,$;\ $P_s(a,b)$ is a polynomial
in $a,b$. Note that $\ (\hh_s)\vert_{q=1}\ =\ (k\cdot x)^s\ $
and thus $\ (\widehat{\exp}_q (k, x))\vert_{q=1}\ =\ \exp (k\cdot
x)\ $. This q-plane wave has some properties analogous
to the classical one but is not an exponent or q-exponent, cf.
 \cite{GR}. This is enabled also by the fact (true also for $q=1$)
that solving the equations may be done in terms of the components
$\ \hh_s\ $. This deformation of the plane wave generalizes
the original one from \cite{DK} to obtain which
one sets $\ P_s(a,b) =0\,$, in which case we shall use
the notation $\ f_s\ $ for the components from \cite{DK} since:
\eqn{fs} \left(\hh_s\right)_{P_s(a,b)=0} ~=~ f_s \ .\ee
Each $\hh_s$ satisfies the q-d'Alembert equation (\ref{qdal})
on the momentum q-cone:
\be\label{zap}\cl^k_q \ \ \equiv \ \ k_- k_+ \ - \ q^{-1}\,
k_v k_{\bv} \ = \ k_+ k_- \ - \ q\, k_v k_{\bv} \ =\ 0 \ .
\ee

The second deformation is given in the basis (\ref{spcb}):
\be\label{copl}
\widetilde{\exp}_q (k, x)\ = \
\sum_{s=0}^\infty\ {1\over [s]_q!}\ \tilde h_s \ ,\ee
\bea\label{tfs}
\tilde h_s \ &=& \ \tb^s\
\sum_{a,b,n}\ {(-1)^{s-a-b}\ q^{n(2a+2b-2n-s)\ +\ a(a-s-1)\ +\
b(s-a-b+1)}\ q^{Q_s(a,b)} \over \G_q(a-n+1)\ \G_q(b-n+1)\
\G_q(s-a-b+n+1)\ [n]_q!}\ \times \nn\\ [2mm]
&& \times\ k_{\bv}^n k_+^{a-n} k_-^{b-n} k_v^{s-a-b+n}\
{\bv}^{s-a-b+n} x_+^{b-n} x_-^{a-n} v^n \ , \\ [2mm]
&&\left(\tb^s \right)^{-1} \ = \
\sum_{p=0}^s \ {q^{(p-s)(p-1)+p} \over [p]_q!\ [s-p]_q!} \ , \nn \eea
where $Q_s(a,b)$ are arbitrary polynomials. If the latter are zero
then $\widetilde{\exp}_q (k, x)$ becomes the q-plane wave
deformation found in \cite{DGPZ}.
The $\th_s$ have the same properties as the $\hh_s$ but the
conjugated basis is used; in particular, they satisfy
the q-d'Alembert equation (\ref{qqdal})
on the momentum q-cone (\ref{zap}).

\section{Solutions of the q-Maxwell equations}
\setcounter{equation}{0}

First we shall use the basis (\ref{spca}).
The solutions of (\ref{mxg}) for $n=0$ in
the homogeneous case ($J=0$) are:
\be\label{resc}
\hF^{h\pm}\ \doteq \ \( _qF^\pm_{0}\)_{J=0}
\ =\ \sum_{m,s=0}^\infty\ {1\over [s]_q!}\
\hF^{h\pm}_{ms} (k) \ f_s\ ,\ee
\bea\label{rescc}
\hF^{h+}_{ms}(k) \ &=&\ \sum_{i=0}^{m}\ \Bigl(\
\sum_{j=0}^{m-i}\
\hp^{ms1}_{ij}\ k_v^i k_-^{m-i-j} k_\bv^j\
(k_v-q^{s+6}zk_-)(k_v-q^{s+3}zk_-) \ + \nn\\[2mm]
&&+\ \hp^{ms2}_{i}\ k_v^i k_\bv^{m-i}\
(k_v-q^{s+6}zk_-)(k_+-q^{s+3}zk_\bv) \ + \nn\\[2mm]
&&+\ \sum_{j=0}^{m-i}\ \hp^{ms3}_{ij}\ k_v^i k_+^{m-i-j}
k_\bv^j\ (k_+-q^{s+6}zk_\bv)(k_+-q^{s+3}zk_\bv)\ \Bigr)\ , \\
\hF^{h-}_{ms}(k) \ &=&\ \sum_{i=0}^{m}\ \Bigl(\
\sum_{j=0}^{m-i}\
\hr^{ms1}_{ij}\ k_v^i k_-^{m-i-j} k_\bv^j\
(k_{\bv}-q^{-1}\bz k_-)(k_\bv -\bz k_-)
\ + \nn\\[2mm] &&+\
\hr^{ms2}_{i}\ k_v^i k_\bv^{m-i}\
(k_+-q^{-1}\bz k_v)(k_{\bv}-\bz k_-)
\ + \nn\\[2mm]
&&+\ \sum_{j=0}^{m-i}\
\hr^{ms3}_{ij}\ k_v^i k_+^{m-i-j} k_\bv^j\
(k_+-q^{-1}\bz k_v)(k_+ -\bz k_v)
\ \Bigr)\label{sashp}\ , \eea
where $\hp^{msa}_{i(j)},\hr^{msa}_{i(j)}$ are independent
constants. The check that these are solutions is done for
commutative Minkowski coordinates and noncommutative momenta on
the q-cone. The terms with $m=0$ of the solutions (\ref{resc}),
(\ref{rescc}), (\ref{sashp}), were obtained earlier \cite{DZ}
(later they were generalized using more general $q$-plane
waves \cite{DPZa}). The solution (\ref{sashp}) can be written
in terms of the deformed plane wave if we suppose that the
$\hr^{msa}_{i(j)}$ for different $s$ coincide:
\ $\hr^{msa}_{i(j)} = \hr^{ma}_{i(j)}\,$. Then we have:
\be\label{tashp} \hF^{h-} \ =\ \sum_{m=0}^\infty\ \hF^{h-}_{m}(k)\
{\exp}_q (k, x) \ , \quad \hF^{h-}_{m}(k) = \hF^{h-}_{ms}(k)\ . \ee

In the inhomogeneous case the solutions of (\ref{mxg}) for $n=0$
are: \bea
_qJ^0 \ &=&\ \bz z \hJ_+ + z \hJ_v + \bz \hJ_\bv + \hJ_- \ ,
\label{current} \\
\hJ_\k \ &=&\
\sum_{m,s=o}^{\infty }{1\over {[s]_q!}}\ \hJ^{ms}_\k(k)\,
f_{s-1} \ ,\quad \k = \pm, v,\bv \ , \label{currena}\\ [2mm]
\hJ^{ms}_+(k)\ &=&\ -\hK^s_{m}(k)\,k_- \ , \label{currea}\\[2mm]
\hJ^{ms}_-(k)\ &=&\ -q^{-s-2}\hK^s_{m}(k)\, k_+ \ , \nn\\[2mm]
\hJ^{ms}_v(k)\ &=&\ \hK^s_{m}(k)\,k_{\bv}\ , \nn\\[2mm]
\hJ^{ms}_{\bv}(k)\ &=&\ q^{-s-2}\hK^s_{m}(k)\,k_v \ , \nn\\[2mm]
\hK^s_{m}(k)\ &\doteq &\ \hg^s_v
k^{m+1}_v+\hg^s_{-}k^{m+1}_-+\hg^s_{+}k^{m+1}_++
\hg^s_{\bv} k^{m+1}_{\bv} \ ; \nn \eea
\bea\label{resf}
_qF^\pm_{0}\ &=& \ \hF^\pm\ +\ \hF^{h\pm}\ , \\
\hF^\pm \ &=&\ \sum_{m,s=0}^\infty\ {1\over [s]_q!}\
\hF^\pm_{ms} (k)\ f_s \ ,\\[2mm]
\hF^+_{ms}(k)\ &=&\ 2d_sq^{-s}\big(
(q^{-s-5}\hg^s_{-}k^m_-+z\hg^s_{v}k^m_v)(k_v-q^{s+3}zk_-)\ + \nn\\[2mm]
&&+\ (q^{-s-5}\hg^s_{\bv}k^m_{\bv}+z\hg^s_{+}k^m_+)(k_+-q^{s+3}z
k_{\bv})\big)\ , \nn\\[2mm]
\hF^-_{ms}(k)\ &=&\ 2d_sq^{-2s-2}\big(
(\hg^s_{-}k^m_-+q^{-2}\bz \hg^s_{\bv}k^m_{\bv})(
k_{\bv }-\bz k_-)\ + \nn\\[2mm]
&&+\ (\hg^s_{v}k^m_v+q^{-2}\bz \hg^s_{+}k^m_+)(k_+-\bz k_v)
\big) \ , \nn \eea
where $d_s = \b^s/\b^{s+1}\ $.
As in the homogeneous case we can not make $\hF^+_{ms}(k)$
independent of $s$. We can make $\hF^-_{ms}(k)$
independent of $s$ by choosing $\hg^s_\k \sim q^{2s}\,d_s^{-1}\,$,
but we can not make ~$\hJ^{ms}_\k(k)$~ independent of $s$.

Since we work with the full Maxwell equations we have also to
check the q-deformation of the current conservation
~$\pd^\nu J_\nu ~=~ 0$~:
\bea  I_{13}\,J ~&=&~ 0 \ ,\label{div} \\
I_{13} ~&=&~ q^3\,[N_z-1]_q\,T_z\,\hd_{\bz}\,\hd_v\,T_v\,T_-\,T_+
~+~ q\,\hd_z\,T_z\,\hd_{\bz}\,\hd_-\,T_v\,T_+ ~+
\nn\\ && +~
q\,[N_z-1]_q\,T_z\,[N_{\bz}-1]_q\,\hd_+\,T_+\,T_{\bv} ~+\nn\\
&&+~ q^{-1}\,[N_{\bz}-1]_q\,\hd_z\,T_z\,\hd_{\bv}\,T_v\,T_-^{-1}\,T_+ ~-\nn\\
&&-~ \l\,
\hm_v\,[N_{\bz}-1]_q\,\hd_z\,T_z\,\hd_-\,\hd_+\,T_v\,T_-^{-1}\, T_+\,T_{\bv}
\label{diva}
\eea
Substituting (\ref{current},\ref{currena}) in the above we get:
\be \label{divb}
q\,J^s_+(k)\,k_+ ~+~ J^s_v(k)\,k_v ~+~
q^{s+2}\,J^s_{\bv}\,k_{\bv} ~+~ q^{s+1}\, J^s_-(k)\,k_-=0
\ee
The latter is fulfilled by the explicit expressions in
(\ref{currea}), but we should note that these expressions
fulfil also the following splittings of (\ref{divb})~:
\bea
&&q\,J^s_+(k)\,k_+ ~+~ J^s_v(k)\,k_v ~=~ 0 \ , \quad
 q\,J^s_{\bv}(k)\,k_{\bv} ~+~ J^s_-(k)\, k_- ~=~ 0 \ , \qquad
\\
&& J^s_+(k)\,k_+ ~+~ q^{s+1}\,J^s_{\bv}(k)\, k_{\bv} ~=~ 0 \ , \quad
J^s_v(k)\,k_v ~+~ q^{s+1}\, J^s_-(k)\, k_- ~=~ 0\ . \qquad
\nn\eea
Furthermore the expressions from (\ref{currea}) fulfil also:
\bea
&&q\,J^s_+(k)\,k_\bv ~+~ J^s_v(k)\,k_- ~=~ 0 \ , \quad
q\,J^s_{\bv}(k)\,k_{+} ~+~ J^s_-(k)\, k_v ~=~ 0 \ , \\
&& J^s_+(k)\,k_v ~+~ q^{s+1}\,J^s_{\bv}(k)\, k_{-} ~=~ 0 \ , \quad
J^s_v(k)\,k_+ ~+~ q^{s+1}\, J^s_-(k)\, k_\bv ~=~ 0\ . \qquad
\nn\eea

Now we shall use the basis (\ref{spcb}).
Then solutions of (\ref{mxg}) for $n=0$ in
the homogeneous case ($J=0$) are:
\be\label{tresc}
\tF^{h\pm}\ \doteq \ \( _qF^\pm_{0}\)_{J=0}
\ =\ \sum_{m,s=0}^\infty\ {1\over [s]_q!}\
\tF^{h\pm}_{ms} (k) \ \th_s\ ,\ee
\bea\label{trescc}
\tF^{h+}_{ms}(k) \ &=&\ \sum_{i=0}^{m}\ \Bigl(\
\sum_{j=0}^{m-i}\
\tp^{ms1}_{ij}\ k_\bv^i k_-^{m-i-j} k_v^j\
(k_v-zk_-)(k_v-qzk_-) \ + \nn\\[2mm]
&&+\
\tp^{ms2}_{i}\ k_\bv^i k_v^{m-i}\ (k_+-zk_\bv) (k_v-qzk_-)
\ + \nn\\[2mm]
&&+\ \sum_{j=0}^{m-i}\
\tp^{ms3}_{ij}\ k_\bv^i k_+^{m-i-j} k_v^j\
(k_+-zk_\bv)(k_+-qzk_\bv)\ \Bigr)\ , \\
\tF^{h-}_{ms}(k) \ &=&\ \sum_{i=0}^{m}\ \Bigl(\
\sum_{j=0}^{m-i}\
\tr^{ms1}_{ij}\ k_\bv^i k_-^{m-i-j} k_v^j\
(k_{\bv}-q^{s+1}\bz k_-)(k_\bv -q^{s+2}\bz k_-)
\ + \nn\\[2mm]
&&+\
\tr^{ms2}_{i}\ k_v^i k_\bv^{m-i}\
(k_{\bv}-q^{s+1}\bz k_-) (k_+-q^{s+2}\bz k_v)
\ + \nn\\[2mm]
&&+\ \sum_{j=0}^{m-i}\
\tr^{ms3}_{ij}\ k_v^i k_+^{m-i-j} k_\bv^j\
(k_+-q^{s+1}\bz k_v)(k_+ -q^{s+2}\bz k_v)
\ \Bigr)\label{tsashp}\  \eea
where $\tp^{msa}_{i(j)},\tr^{msa}_{i(j)}$ are independent
constants, ~$Q_s(a,b)=0$~ in ~$\th_s\,$.
The terms with $m=0$ of the solutions (\ref{tresc}), (\ref{trescc}),
(\ref{tsashp}) were obtained earlier in \cite{DZ}
(and using the generalized q-plane wave in \cite{DPZa}).
The solution (\ref{trescc}) can be written
in terms of the deformed plane wave if we suppose that
the $\tp^{msa}_{i(j)}$ for different $s$ coincide:
\ $\tp^{msa}_{i(j)} = \tp^{ma}_{i(j)}\,$. Then we have:
\be\label{ttashp} \tF^{h+} \ =\ \sum_{m=0}^\infty\ \tF^{h+}_{m}(k)\
\widetilde{\exp}_q (k, x) \ , \quad \tF^{h+}_{m}(k) =
\tF^{h+}_{ms}(k)\ . \ee

In the inhomogeneous case the solutions of (\ref{mxg}) for $n=0$
are: \bea
_qJ^0 \ &=&\ \bz z \tJ_+ + z \tJ_v + \bz \tJ_\bv + \tJ_- \ ,
\label{tcurrent} \\
\tJ_\k \ &=&\
\sum_{m,s=o}^{\infty }{1\over {[s]_q!}}\ \tJ^{ms}_\k(k)\,
\th_{s-1} \ ,\quad \k = \pm, v,\bv \ ,\label{tcurrena}\\ [2mm]
\tJ^{ms}_+(k)\ &=&\ -q^{s+1}\tK^s_{m}(k)\,k_- \ , \label{tcurrea}
\\[2mm]
\tJ^{ms}_-(k)\ &=&\ -q^{-1}\tK^s_{m}(k)\, k_+ \ , \nn\\[2mm]
\tJ^{ms}_v(k)\ &=&\ \tK^s_{m}(k)\,k_{\bv}\ , \nn\\[2mm]
\tJ^{ms}_{\bv}(k)\ &=&\ q^{s}\tK^s_{m}(k)\,k_v \ , \nn\\[2mm]
\tK^s_{m}(k)\ &\doteq &\ \tg^s_v
k^{m+1}_v+\tg^s_{-}k^{m+1}_-+\tg^s_{+}k^{m+1}_++ \tg^s_{\bv}
k^{m+1}_{\bv} \ , \nn\eea
\bea\label{tresf}
_qF^\pm_{0}\ &=& \ \tF^\pm\ +\ \tF^{h\pm}\ , \\
\tF^\pm \ &=&\ \sum_{m,s=0}^\infty\ {1\over [s]_q!}\
\tF^\pm_{ms} (k)\ \th_s \ ,\\[2mm]
\tF^+_{ms}(k)\ &=&\ 2\td_sq^{s-2}\big(
(\tg^s_{-}k^m_- + q^{-1} z\tg^s_{v}k^m_v)(k_v-qzk_-)\ + \nn\\[2mm]
&&+\ (\tg^s_{\bv}k^m_{\bv}+ q^{-1}z\tg^s_{+}k^m_+)(k_+-qz
k_{\bv})\big)\ , \nn\\[2mm]
\tF^-_{ms}(k)\ &=&\ 2\td_s \big(
(q^{-s-3} \tg^s_{-}k^m_- + q \bz \tg^s_{\bv}k^m_{\bv})(
k_{\bv }-q^{s+2}\bz k_-)\ + \nn\\[2mm]
&&+\ (q^{-s-3} \tg^s_{v}k^m_v + q \bz \tg^s_{+}k^m_+)(k_+-q^{s+2} \bz k_v)
\big) \ , \nn \eea
where $\td_s = \tb^s/\tb^{s+1}\,$, ~$Q_s(a,b)=0$~ in ~$\th_s\,$.
We can not make ~$\tF^-_{ms}(k)$~ or ~$\tJ^{ms}_\k(k)$~
independent of $s$. We can make $\tF^+_{ms}(k)$ independent of
$s$ by choosing $\tg^s_\k \sim q^{-s}\,\td_s^{-1}\,$.

Also here we shall check whether  the q-deformation of the
current conservation (\ref{div}) is fulfilled. The analog of
(\ref{diva}) in the basis (\ref{spcb}) is:
\bea
I_{13} ~&=&~
[N_z-1]_q\,\hd_\bz\,T_{\bz}\,\hd_v\,T_\bv\,T_+\,T_-^{-1}
~+~ q\,\hd_\bz\,T_\bz\,\hd_{z}\,\hd_-\,T_\bv\,T_+ ~+
\nn\\ && +~
q\,[N_\bz-1]_q\,T_\bz\,[N_{z}-1]_q\,\hd_+\,T_+\,T_{v} ~+\nn\\
&&+~
q^{2}\,[N_{\bz}-1]_q\,\hd_z\,T_\bz\,\hd_{\bv}\,T_\bv\,T_-\,T_+
~-\nn\\ &&-~ \l\,q\,
\hm_v\,[N_{\bz}-1]_q\,\hd_z\,T_\bz\,\hd_-\,\hd_+\,T_-\, T_+
\label{divaa}
\eea
Then the analog of (\ref{divb}) is:
\be \label{divc}
J^s_+(k)\,k_+ ~+~ q^s\, J^s_v(k)\,k_v ~+~
J^s_{\bv}\,k_{\bv} ~+~ q^{s}\, J^s_-(k)\,k_-=0
\ee
The latter is fulfilled by the explicit expressions in
(\ref{tcurrea}), but we should note that these expressions
fulfil also the following splittings of (\ref{divc})~:
\bea
&&J^s_+(k)\,k_+ ~+~ q^s\, J^s_v(k)\,k_v ~=~ 0 \ , \quad
J^s_{\bv}(k)\,k_{\bv} ~+~ q^s\, J^s_-(k)\, k_- ~=~ 0\ , \qquad
\nn\\
&& J^s_+(k)\,k_+ ~+~ J^s_{\bv}(k)\, k_{\bv} ~=~ 0 \ , \quad
J^s_v(k)\,k_v ~+~  J^s_-(k)\, k_- ~=~ 0\ . \qquad
\eea
Furthermore the expressions from (\ref{tcurrea}) fulfil also:
\bea
&&J^s_+(k)\,k_\bv ~+~ q^s\, J^s_v(k)\,k_- ~=~ 0 \ , \quad
  J^s_{\bv}(k)\,k_{+} ~+~ q^s\, J^s_-(k)\, k_v ~=~ 0\ , \qquad \nn\\
&& J^s_+(k)\,k_v ~+~ J^s_{\bv}(k)\, k_{-} ~=~ 0 \ , \quad
J^s_v(k)\,k_+ ~+~  J^s_-(k)\, k_\bv ~=~ 0 \ .\qquad
\eea

Summarizing,
we have given new solutions of the full q-Maxwell equations in
two conjugated bases (\ref{spca}) and (\ref{spcb}). The solutions
of the homogeneous equations are also new (the old solutions are
special cases).  We see that the roles of the solutions $F^{+}$
and $F^{-}$ are exchanged in the two conjugated bases.
We note also that the currents components  are different~: ~$\hJ^{ms}_\k
~\neq~ \tJ^{ms}_\k$~ (for $q\neq 1$, $\k\neq v$), and in both cases they
can not be made independent of $s$. Thus, there is no advantage
of choosing either of the bases (\ref{spca}) or  (\ref{spcb}). It
may be also possible to use both in a Connes-Lott type model
\cite{ConLott}.

\section{Linear conformal gravity}
\setcounter{equation}{0}

We consider now the quantum group  analogs of linear conformal
gravity following the approach of \cite{VKD4}.
We start with the $q=1$ situation and we
first write the Weyl gravity equations in an indexless formulation,
trading the indices for two conjugate variables $z,\bz$,
just as  for the Maxwell equations.

Weyl gravity is governed by the Weyl tensor:
\eqn{wt} C_{\mu\nu\s\t} = R_{\mu\nu\s\t} -\half(
g_{\mu\s}R_{\nu\t}+
g_{\nu\t}R_{\mu\s}-
g_{\mu\t}R_{\nu\s}-
g_{\nu\s}R_{\mu\t}) +\sixth (g_{\mu\s}g_{\nu\t}-
g_{\mu\t}g_{\nu\s})R \ ,\ee
where $g_{\mu\nu}$ is the metric tensor.
Linear conformal gravity is obtained when the metric tensor
is written as: $g_{\mu\nu} = \y_{\mu\nu} + h_{\mu\nu}$,
where $\y_{\mu\nu}$ is the flat Minkowski metric, $h_{\mu\nu}$
are small so that all quadratic and higher order terms are neglected.
In particular:  $R_{\mu\nu\s\t} = \half(\pd_\mu \pd_\t h_{\nu\s}+
\pd_\nu \pd_\s h_{\mu\t}-
\pd_\mu \pd_\s h_{\nu\t}-
\pd_\nu \pd_\t h_{\mu\s})$.
The equations of linear conformal gravity are:
\eqn{wfr}  \pd^\nu \pd^\t C_{\mu\nu\s\t} = T_{\mu\s}\ , \ee
where  $T_{\mu\nu}$ is the energy-momentum tensor.
{}From the symmetry properties of the Weyl tensor it follows that
it has ten independent components. These may be chosen as follows
(introducing notation for future use):
\eqnn{weyli}
&C_0=C_{0123}\ , \quad
C_1=C_{2121}\ ,\quad
C_2=C_{0202}\ ,\quad
C_3=C_{3012}\ ,\nn\\
&C_4=C_{2021}\ ,\quad
C_5=C_{1012}\ ,\quad
C_6=C_{2023}\ ,\nn\\
&C_7=C_{3132}\ ,\quad
C_8=C_{2123}\ ,\quad
C_9=C_{1213}\ . \eea
Furthermore, the Weyl tensor  transforms
as the direct sum of two conjugate Lorentz irreps, which
we shall denote as $C^\pm$.
The tensors $T_{\mu\nu}$ and $h_{\mu\nu}$ are symmetric and traceless
with nine independent components.

In order to be more precise we recall that the physically relevant
representations $T^\chi$ of the
4-dimensional conformal algebra $su(2,2)$ may be labelled by
$\chi = [n_1,n_2;d]$, where $n_1, n_2$ are non-negative
integers fixing finite-dimensional irreducible representations of
the Lorentz subalgebra, (the dimension being
$(n_1 +1)(n_2 +1)$), and ~$d$~ is the conformal
dimension (or energy). (In the literature these Lorentz
representations are labelled also by $(j_1,j_2)=(n_1/2,n_2/2)$.)
The Weyl tensor transforms as the direct sum:
\eqnn{sgw} &\chi^+ \oplus \chi^- \cr &\cr
&\chi^+ = [4,0;2] ~,
\quad    \chi^- = [0,4;2] ~, \eea
while the energy-momentum tensor and the metric transform as:
\eqn{sgt} \chi_T = [2,2;4] ~, \qquad
\chi_h = [2,2;0] ~, \ee
as anticipated. Indeed, $(n_1,n_2) = (2,2)$ is the nine-dimensional
Lorentz representation, (carried by $T_{\mu\nu}$ or $h_{\mu\nu}$), and
$(n_1,n_2) = (4,0),(0,4)$ are the two conjugate five-dimensional
Lorentz representations, (carried by $C^\pm$),
while the conformal dimensions are the canonical dimensions
of a energy-momentum tensor ($d=4$), of the metric ($d=0$),
and of the Weyl tensor ($d=2$). (For comparison, note
that the Maxwell components
$F^+$, $F^-$, used in the previous sections, have signatures:
$[2,0;2]$, $[0,2;2]$, resp., while the current $J$ has
signature $[1,1;3]$.)
Further, we shall use again the fact that a Lorentz irrep (spin-tensor)
with signature $(n_1,n_2)$ may be represented by
a polynomial $G(z,\bz)$ in $z,\bz$ of order $n_1,n_2$, resp.
More explicitly, for the irreps mentioned above we use:
\eqna{deff}
C^+(z) &=& z^4C^+_4+z^3C^+_3+z^2C^+_2+zC^+_1+C^+_0 \ ,\\
C^-(\bz) &=& \bz^4C^-_4+\bz^3C^-_3+\bz^2C^-_2+\bz C^-_1+C^-_0 \ ,\\
T(z,\bz )& =& z^2\bz^2T'_{22}+z^2\bz T'_{21}+z^2T'_{20}+\cr
&& +z\bz^2T'_{12}+z\bz T'_{11}+zT'_{10}+ \cr
&& +\bz^2T'_{02}+\bz T'_{01}+T'_{00} \ ,\\
h(z,\bz ) &=& z^2\bz^2h'_{22}+z^2\bz h'_{21}+z^2h'_{20}+\cr
&& +z\bz^2h'_{12}+z\bz h'_{11}+zh'_{10}+\cr
&& +\bz^2h'_{02}+\bz h'_{01}+h'_{00} \ , \eena
where the indices on the RHS are not Lorentz-covariance indices,
they just indicate the powers of $z,\bz$. The components
$C^\pm_k$ are given in terms of the Weyl tensor components as follows:
\eqnn{wcomp}
&&C^+_0=C_2-\ha C_1-C_6+i(C_0+ \ha C_3+
C_7)\cr
&&C^+_1=2(C_4-C_8+i(C_9-C_5))\cr
&&C^+_2=3(C_1-iC_3)\cr
&&C^+_3=8(C_4+C_8+i(C_9+C_5))\cr
&&C^+_4=C_2- \ha C_1+C_6+i(C_0+ \ha C_3 -
C_7)\cr
&&C^-_0=C_2-{C_1\over 2}-C_6-i(C_0+ \ha C_3+
C_7)\cr
&&C^-_1=2(C_4-C_8-i(C_9-C_5))\cr
&&C^-_2=3(C_1+iC_3)\cr
&&C^-_3=2(C_4+C_8-i(C_9+C_5))\cr
&&C^-_4=C_2- \ha C_1 +C_6-i(C_0+ \ha C_3-
C_7) \eea
while the components $T'_{ij}$ are given in terms of $T_{\mu\nu}$
as follows:
\eqnn{enem}
&&T'_{22} = T_{00} + 2 T_{03} + T_{33} \cr
&&T'_{11} = T_{00} - T_{33} \cr
&&T'_{00} = T_{00} - 2 T_{03} + T_{33} \cr
&&T'_{21} = T_{01} + iT_{02} + T_{13} + iT_{23} \cr
&&T'_{12} = T_{01} - iT_{02} + T_{13} - iT_{23} \cr
&&T'_{10} = T_{01} + iT_{02} - T_{13} - iT_{23} \cr
&&T'_{01} = T_{01} - iT_{02} - T_{13} + iT_{23} \cr
&&T'_{20} = T_{11} + 2iT_{12} - T_{22} \cr
&&T'_{02} = T_{11} - 2iT_{12} - T_{22} \eea
and similarly for $h'_{ij}$ in terms of $h_{\mu\nu}\,$.

In these terms all linear conformal (Weyl) gravity equations \eqref{wfr} may be written in compact form as the following pair of  equations:
\eqna{wge}
\tI^+ ~ C^+(z)
~=&~ T(z,\bz) ~, \\
\tI^- ~ C^-(\bz)
~=&~ T(z,\bz) ~, \eena
where  the operators $I^\pm$ are given as follows:
\eqna{oper}
\tI^+ &=&
\Bigl(z^2\bz^2\pd_+^2+z^2\pd^2_v+\bz^2\pd_{\bv }^2+\pd^2_-+\cr
&& +2z^2\bz \pd_v\pd_++2z\bz^2\pd_+\pd_{\bv }+2z\bz
(\pd_-\pd_++\pd_v\pd_{\bv })+\cr
&& +2\bz \pd_-\pd_{\bv }+2z\pd_v\pd_-\Bigr)\pd^2_z -\cr
&& -6\Bigl(z\bz^2\pd_+^2+z\pd^2_v+2z\bz \pd_v\pd_++\bz^2\pd_+\pd_{\bv }+\cr
&& +\bz(\pd_-\pd_++\pd_v\pd_{\bv })+\pd_v\pd_-\Bigr)\pd_z+\cr
&& 12\Bigl(\bz^2\pd_+^2+\pd_v^2+2\bz\pd_v\pd_+\Bigr)\ , \\
\tI^{-} &=&
\Bigl(z^2\bz^2\pd_+^2+z^2\pd^2_v+\bz^2\pd_{\bv }^2+\pd^2_-+\cr
&& +2z^2\bz \pd_v\pd_++2z\bz^2\pd_+\pd_{\bv }+2z\bz
(\pd_-\pd_++\pd_v\pd_{\bv })+\cr
&& +2\bz \pd_-\pd_{\bv }+2z\pd_v\pd_-\Bigr)\pd^2_{\bz } -\cr
&& -6\Bigl(z^2\bz \pd_+^2+\bz\pd^2_{\bv }+2z\bz \pd_+\pd_{\bv}+
z^2\pd_v\pd_++\cr
&& +z(\pd_-\pd_++\pd_v\pd_{\bv })+\pd_-\pd_{\bv }\Bigr)\pd_{\bz }+ \cr
&& 12\Bigl(z^2\pd_+^2+\pd_{\bv }^2+2z\pd_+\pd_{\bv }\Bigr)
\ . \eena

To make more transparent the origin of these expressions and in
the same time to derive the quantum group deformation of \eqref{wge},
\eqref{oper} we first introduce the following parameter-dependent operators:
\eqna{opsb}
\tI^+(n) ~&=&~ \ha \Bigl( n(n-1) I^2_1 I^2_2
- 2(n^2-1) I_1 I^2_2 I_1  + n(n+1) I^2_2 I^2_1
\Bigr)
~, \\
\tI^-(n) ~&=&~ \ha \Bigl( n(n-1) I^2_3 I^2_2
- 2(n^2-1) I_3 I^2_2 I_3  + n(n+1) I^2_2 I^2_3
~, \eena
where
\eqn{opsbb}
 I_1 ~\equiv ~\pd_z ~, \quad I_2 ~\equiv ~ \bz z \pd_+ + z\pd_v
+ \bz \pd_\bv + \pd_-  ~, \quad I_3 ~\equiv ~\pd_\bz ~. \ee
It is easy to check that we have the following relation:
\eqn{rel} \tI^\pm = \tI^\pm(4) \ . \ee
We note in passing that group-theoretically the operators $I_a$
correspond to the three simple roots of the root system of
$sl(4)$, while the operators $I^\pm_n$ correspond to the two
non-simple non-highest roots \cite{VKDcl}.

This is the form that is immediately generalizable to the
$q$-deformed case. Using results from  \cite{VKD4} we have:
\eqna{opsc} _q\tI^+(n) ~&=&~ \ha \Bigl( [n]_q\,[n-1]_q\ _qI^2_1
\,_qI^2_2 - [2]_q\,[n-1]_q\,[n+1]_q \,_qI_1 \,_qI^2_2 \,_qI_1 +\cr
&&+ [n]_q\,[n+1]_q \,_qI^2_2 \,_qI^2_1 \Bigr)
~, \\
_q\tI^-(n) ~&=&~ \ha \Bigl( [n]_q\,[n-1]_q \,_qI^2_3 \,_qI^2_2 -
[2]_q\,[n-1]_q\,[n+1]_q \,_qI_3 \,_qI^2_2 \,_qI_3  + \cr&&+
[n]_q\,[n+1]_q \,_qI^2_2 \,_qI^2_3 \Bigr) ~, \eena where the
$q$-deformed versions $\,_qI_a$ of \eqref{opsbb} in the basis
(\ref{spca})  are: \eqna{rek}
_qI_1\ &=&\ \hd_zT_zT_vT_+(T_-T_{\bv})^{-1} \\
_qI_2\ &=&\ (q\hm_z\hd_vT_-^2+\hm_z\hm_{\bz}\hd_+T_-T_{\bv}T_v^{-1}+
\hd_-T_-\ +\cr
&&+\ q^{-1}\hm_{\bz}\hd_{\bv}-\l\hm_v\hm_{\bz}\hd_-\hd_+T_{\bv})\ T_{\bv}
T_{\bz}^{-1} \\
_qI_3\ &=&\ \hd_{\bz}T_{\bz}\ . \eena Then the $q$-Weyl equations
are: \eqna{wgq} _q\tI^+(4) ~ C^+(z)
~=&~ T(z,\bz) ~, \\
_q\tI^-(4) ~ C^-(\bz)
~=&~ T(z,\bz) ~. \eena
(For comparison, note that for the derivation of
the $q$-Maxwell operators \eqref{mxg} were used the
following expressions:
$\,_qI^+_n = \ha ( [n+2]_q\, _qI_1\,  _qI_2  -
[n+3]_q\, _qI_2\, _qI_1 )$,
$\,_qI^-_n ~=~ \ha ( [n+2]_q\, _qI_3\, _qI_2  -
[n+3]_q\, _qI_2\, _qI_3 )$.)

Finally, we write down the pair of equations which give the Weyl tensor
components in terms of the metric tensor:
\eqna{wme}
_q\tI^+(2) ~ h(z,\bz)
~=&~ C^+(z) ~, \\
_q\tI^-(2) ~ h(z,\bz)
~=&~ C^-(\bz)
~. \eena

We stress the advantage of the indexless formalism due to which
two different pairs of equations, \eqref{wgq},
\eqref{wme}, may be written using the same parameter-dependent
operator expressions by just specializing the values of
the parameter.

\section*{Acknowledgments}
The authors would like to thank the Alexander von Humboldt
Foundation for financial support in the framework of the
Clausthal-Leipzig-Sofia Cooperation,  and the Bulgarian National Council
for Scientific Research for financial support under
grant F-1205/02.

\end{document}